\documentclass[12pt]{amsproc}
\usepackage{amssymb}
\usepackage[mathscr]{eucal}
\usepackage{amscd}

\usepackage[dvips]{color, graphics}

\usepackage[cp850]{inputenc}
\usepackage{amssymb}
\usepackage[mathscr]{eucal}
\usepackage{wrapfig}

\title{Seven views on approximate convexity and
the geometry of $K$-spaces} 
\thanks{The research of the first and second named authors is 
supported in part by
DGICYT---project-BMF-2001-0813}
\author{F\'elix Cabello S\'anchez, Jes\'us M. F. Castillo and~Pier~L.~Papini}
\address{Departamento de Matem\'aticas, Universidad de Extremadura,
06071-Badajoz, Spain}
\email{fcabello@unex.es, castillo@unex.es}
\address{Dipartimento di Matematica, Universita di Bologna,
Piazza Porta S. Donato 5, Bologna, Italy}
\email{papini@dm.unibo.it}

\begin{document}
\noindent{\footnotesize Version: April 2005}
\bigskip

\bigskip

\maketitle
\newtheorem{proposition}{Proposition}
\newtheorem{definition}{Definition}
\newtheorem{corollary}{Corollary}
\newtheorem{question}{Question}
\newcommand{\R}{\mathbb{R}}
\newcommand{\e}{\varepsilon}
\newcommand{\N} {\mathbb{N}}
\newcommand{\vr}{\varrho}
\newtheorem{theorem}{Theorem}
\newtheorem{lemma}{Lemma}
\newtheorem{example}{Example}
\newtheorem{remark}{Remark}

\def\co{\operatorname{co}}
\def\sgn{\operatorname{sgn}}
\def\dist{\operatorname{dist}}
\def\epi{\operatorname{epi}}
\def\supp{\operatorname{supp}}

\bibliographystyle{plain}

\begin{abstract}
We study the interplay between the behaviour of
approximately convex (and approximately affine) functions on 
the unit ball of a Banach space
and the geometry of Banach
$K$-spaces.
\end{abstract}

\section*{Introduction}
This paper deals with the local stability of convexity, affinity
and Jensen functional equation on infinite dimensional Banach
spaces. Recall that a function $f:D\to\R$ is said to be
$\e$-convex if it satisfies
$$
f(tx+(1-t)y)\leq tf(x)+(1-t)f(y)+\e
$$
for all $x,y\in D,t\in[0,1]$. If no specific $\e$ is required we
speak of an approximately convex function. Of course, any
arbitrary function which is uniformly close to a true convex
function is approximately convex. These will be called trivial or
approximable. It may happen that there are no more: Hyers and Ulam
proved in \cite{hyerulam}  that if $D$ is a convex set in $\R^n$,
then for every $\e$-convex function $f:D\to\R$ there exists a
convex $a:D\to\R$ such that
$$
\sup_{x\in D}|f(x)-a(x)|\stackrel{\text{def}}{=} d_D(f,a)\leq
C\cdot\e,
$$
where $C=C_n$ is a constant depending only on $n$. It is
apparent that the papers \cite{hyerulam, gree, chol, cho,
casipapi, dilw, dilw2, dilw3, dilw4, cabe, lacz} contain the
complete story of $C_n$.

As far as we know, the first connections between approximately
convex functions and the geometry of infinite dimensional Banach
spaces appear in \cite{casipapi} and \cite{cabe}. In
\cite{casipapi} it was proved that Lipschitz $\e$-convex functions
are approximable on bounded sets of B-convex spaces, with the distance to the
approximating convex function depending only on $\e$. (Recall that
B-convexity means `having non-trivial type $p>1$'.) That paper
contains some counterexamples based on the fact that $\ell_1$ is
the Banach envelope of the spaces $\ell_p$ for all $0<p<1$. In
\cite{cabe} it was remarked that every Banach space which is not a
$K$-space (see Section \ref{K} for precise definitions) admits a
`bad' (that is, non-approximable) $\e$-convex function defined on
its unit ball. To be a $K$-space is a homological property of
Banach spaces which is closely related to the behaviour of
quasi-linear maps. It suffices to recall here that the spaces
$\ell_p$ and $L_p$ are $K$-spaces if and only if $p\neq 1$.

With this background, let us explain the contents and the
organization of the paper.

The first Section is preliminary: we use the fact that $\ell_1$ is
not a $K$-space to obtain explicit examples of `bad' approximately
convex functions on infinite dimensional simplexes (the examples
in \cite{cabe, dilw2} are not explicit.) This leads to the
question of whether $K$-spaces admit `bad' $\e$-convex functions
on their unit balls. The (affirmative) answer comes in Section 2,
where we exhibit a non-trivial approximately convex function on
the infinite dimensional cube (the unit ball of $\ell_\infty$).
This solves the main problem raised in \cite{cabe}.

Having seen that the local stability of convexity does not hold in
$K$-spaces, we prove in Section 4 that the local stability of
affinity is equivalent to being a $K$-space: precisely, a Banach
space $X$ is a $K$-space if and only if for every $\e$-affine
function $f$ defined on its unit ball $B_X$ there exists a true
affine $a:B_X\to\R$ such that
$$
d_{B_X}(f,a)\leq A\cdot\e,
$$
where $A$ is a constant depending only on $X$. In Section~4 we prove
a similar result for Jensen's functional equation
$$ f\left(\frac{x+y}{2}\right)=\frac{f(x)+f(y)}{2}.
$$
We conclude Sections 3 and 4 by showing that there is a universal
constant $A_D\leq 224$ such that if $f$ is an $\e$-affine function
defined on the $n$-dimensional (euclidean) ball $D$ then there
exists an affine function $a:D\to\R$ with $d_D(f,a)\leq A_D$. A
similar result is proved for $n$-cubes. This solves a problem
posed by Laczkovich in \cite{lacz}.

Section 5 deals with the question of the uniform approximation. Under rather
mild assumptions on the convex set $D$ we show that if every $\e$-affine
(respectively, $\e$-Jensen) $f:D\to\R$ is approximable by an affine
(respectively, Jensen) function $a$, then this can be done with
$d_D(f,a)\leq M\cdot\e$, where $M$ is a constant depending only on $D$.

Finally, Sections 6 and 7 deal with Banach envelopes. In some
sense, the Banach envelope $\co X$ is the nearest Banach space to
a given quasi-Banach space $X$. Here, we regard $X$ as a
topological vector space whose topology is `approximately convex'
while that of $\co X$ is truly convex. Our contribution
complements previous results by Kalton. Precisely, we show that
$c_0$ is not isometric to the Banach envelope of a non-locally
convex space (with separating dual) although there are non-locally
convex spaces $X$ whose Banach envelopes are arbitrarily close (in
the Banach-Mazur distance) to $c_0$ .

\section{Quasi-linear and approximately convex maps}\label{K}

Let $X$ and $Y$ be (real) Banach spaces. A map $f:X\to Y$ is said
to be quasi-linear if:
\begin{itemize}
\item it is homogeneous: $f(tx)=tf(x)$ for all $x\in X$ and
$t\in\R$.
\item it is quasi-additive: $\|f(x+y)-f(x)-f(y)\|\leq
Q(\|x\|+\|y\|)$ for some constant $Q$ independent on $x,y\in X$.
\end{itemize}

The least possible constant in the above inequality is denoted
$Q(f)$ and referred to as the quasi-additivity constant of the map
$f$. When $f$ is homogeneous we also speak of $Q(f)$ as the
quasi-linearity constant of $f$.

Although the original notion of a $K$-space refers to the
possibility of lifting operators (see \cite{kaltpeckrobe} for
background), it will be convenient for our purposes to give the
following definition.

\begin{definition}
\emph{A Banach space is a $K$-space if for every quasi-linear map
$f:X\to\R$ there is a linear (although not necessarily
continuous!) functional $\ell:X\to\R$ such that
$$|f(x)-\ell(x)|\leq M\|x\|$$ for some $M$ and all $x\in X$.}
\end{definition}

Thus $K$-spaces are closely related to the stability of linear
functionals. But that stability is ``asymptotic'' ra\-ther than
``epsilonic''. It will be convenient to introduce the following
asymptotic distance for functions acting between Banach spaces: $$
\dist(f,g)=\inf\{M:\|f(x)-g(x)\|\leq M\|x\|\text{ for all } x\},
$$ where the infimum of the empty set is treated as infinity. In
this way, a $K$-space is a Banach space in which every
quasi-linear functional $f$ is at finite distance from some linear
functional $\ell$. In fact $\ell$ can be chosen in such a way that
$\dist(f,\ell)\leq \kappa\cdot Q(f)$, where $\kappa$ is a constant
depending only on $X$; see \cite[proposition 3.3]{kaltf}.

There are two main types of $K$-spaces: B-convex spa\-ces
\cite{kaltf, kaltpeckrobe, kalttype} and $\mathcal
L_\infty$-spaces \cite{kaltrobe}. Thus, for instance, the
classical spaces $\ell_p$ and $L_p$ are $K$-spaces for all
$1<p\leq\infty$ as well as $c_0$ and all $C(K)$ spaces.

On the other hand, $\ell_1$ (and also every infinite dimensional
$\mathcal L_1$-space)
is not a $K$-space. This was proved by Kalton \cite{kaltf,
kaltstud}, Ribe \cite{ribe} and Roberts \cite{robe}. In fact
Kalton \cite{kaltf} and Ribe \cite{ribe} give (more or less)
explicit examples of quasi-linear maps $f:\ell\to\R$  with
$\dist(f,\ell)=\infty$ for all linear maps $\ell:\ell_1\to\R$.

Ribe's map is given by $$ R(x)=\sum_i
x_i\log_2|x_i|-\left(\sum_ix_i\right)\log_2\left|\sum_ix_i\right|, $$
where $x=\sum_ix_ie_i$ and assuming $0\log 0=0$. It is quasi-linear with
constant 2. Kalton's map is
defined as $$ K(x)=\sum_i\tilde x_i\log i\quad\quad(x\geq0), $$
where $\tilde x$ is the decreasing arrangement of $x$ and then
extended to the finitely supported sequences of $\ell_1$ by $$
K(x)=K(x^+)-K(x^-). $$ It should be noted that the above
formul\ae\ have sense only for finitely supported sequences.
Nevertheless, quasi-linear maps can be extended from dense
subspaces to the whole space (preserving quasi-linearity).

One of the basic observations in \cite{cabe} was that the
restriction of a quasi-linear map to a bounded set of a Banach
space is an approximately convex function which is uniformly close
to a convex function if and only if the starting quasi-linear map
is asymptotically close to a linear map. Since in this paper we
shall deal with convexity, affinity and Jensen functional
equation, let us state  a slightly stronger result for midconvex
functions. Recall that a midconvex function is one satisfying the
inequality $$ f\left(\frac{x+y}{2}\right)\leq\frac{f(x)+f(y)}{2}.
$$

\begin{lemma}\label{midconvex}
Let $f:X\to\R$ be a quasi-linear function. Suppose there is a
midconvex function $a:B_X\to \R$ such that $d_{B_X}(f,a)<\infty$.
Then there is a linear map $\ell:X\to\R$ such that
$\dist(f,\ell)<\infty$.
\end{lemma}

\begin{proof}
This follows from \cite[proof of theorem 2]{cabe} (in which the
result was proved for convex $a$) and \cite[lemma 8.8]{hir}
(asserting that the midconvex function $a$ is actually convex).
\end{proof}

In this way, every (non-trivial) quasi-linear map produces a
``bad'' approximately convex function on the ball of the
corresponding Banach space. Even if nobody knows the values of
Ribe's function on points $x$ for which the series $\sum_n
x_n\log_2|x_n|$ diverges, we can use it to produce an explicit
counterexample on the ``infinite-dimensional simplex'' $$
\Delta^\infty=\left\{x\in\ell_1:x_i\geq0\text{ for all }i \text{
and } \sum_{i=1}^\infty x_i=1\right\}. $$ Indeed, if
$\Delta^\infty_{00}$ denotes the subset of all finitely supported
sequences in $\Delta^\infty$, then $\Delta^\infty_{00}$ is convex
and $\Delta^\infty\backslash \Delta^\infty_{00}$ acts as an ideal:
if $x\in \Delta^\infty\backslash \Delta^\infty_{00}$ and
$y\in\Delta^\infty$ then every non-trivial convex combination
$tx+(1-t)y$ belongs to $\Delta^\infty\backslash
\Delta^\infty_{00}$. Therefore the function
\begin{equation}\label{ribe}
f(x)=\begin{cases} -\sum_{i=1}^\infty x_i\log_2x_i & x\in
\Delta^\infty_{00}\\ 0 & x\notin \Delta^\infty_{00}\\
\end{cases}
\end{equation}
is $1$-convex on $\Delta^\infty$, but
$d_{\Delta^\infty}(f,g)=\infty$ for any convex $g$. Kalton's
function leads to another explicit counterexample taking $ f(x)=
\sum_{i=1}^\infty \tilde x_i\log_2 i$ for $x\in
\Delta^\infty_{00}$ and $f(x)=0$ for $ x\notin
\Delta^\infty_{00}$.

This yields explicit counterexamples for the Hyers-Ulam stability
of convexity in infinite dimensions:

\begin{proposition}[Compare to \cite{cabe} and \cite{dilw}]\label{infinite}
Every infinite dimensional Banach space contains
 a compact convex set $D\subset X$ and
an approximately convex map $h:D\to\R$ such that $d_D(h,g)=\infty$
for every convex $g$.
\end{proposition}

\begin{proof}
Let $(e_n)_n$ be a (normalized) basic sequence in $X$. The map
$\Phi:\Delta^\infty\to X$ sending
$(t_n)\in\Delta^\infty$ into $\sum_{n=1}^\infty t_n(e_n/n)$
defines an one-to-one affine map between $\Delta^\infty$ and a certain
compact convex set $D\subset X$.
Now, if $f$ is the function given in (\ref{ribe}), then $h=f\circ\Phi^{-1}$
is a non-approximable 1-convex function on $D$.
\end{proof}

Let $f:D\to\R$ be an arbitrary function defined on a convex set. The function
$$
\co f (x)= \inf\left\{\sum_{i=1}^nt_if(x_i):
x=\sum_{i=1}^nt_ix_i\right\}
$$
represents the greatest convex minorant of $f$. It is clear that $f$ is
uniformly close to some convex function on $D$ if and only if ($\co f$ takes
only finite values and) $d_D(f,\co
f)<\infty$. Actually the distance from $f$ to the convex functions is
$\frac{1}{2}d_D(f,\co f)$ and  it is attained at $g=\co
f+\frac{1}{2}d_D(f,\co f)$. This will be used without further mention.

\section{Bad $\e$-convex functions on the unit ball
of $\ell_\infty$}

All non-trivial $\e$-convex functions presented so far depend on
the fact that $\ell_1$ is not a $K$-space. So, at this moment the
question is if there are bad $\e$-convex functions on the ball of
a $K$-space. To tackle this question one needs a completely
different type of approximately convex function: since we will
prove in Section \ref{sec:af} that the notion of affinity is
stable on the ball of a $K$-space, one idea is to work with a
function which is approximately convex, but not approximately
concave.

Such an example is supplied by Cholewa and Kominek \cite{cho,ko}
(see also \cite{bz}) as follows: let $c^+_{00}$ be the positive
cone of the space of all finitely supported sequences; for $x\in
c^+_{00}$ set $m(x) = \max_i  x_i $ and then $$ \omega(x) =
\mathrm{min} \lbrace n \in \N : m(x) \geq 2^{-n} \rbrace.$$
This
function is $2$-convex on $c^+_{00}$ (\cite{cho,ko}). Since
$\omega(e_n)=0$ but $$\omega \left( \frac{1}{2^n} \sum_{i=1}
^{2^n}
 e_i \right) = n,$$
the function $\omega$ is not uniformly close to any convex
 function on $c_{00}^+$.
 Since every infinite dimensional space (no topology is assumed)
  contains a subset affinely isomorphic to $c_{00}^+$ we obtain a
 non-trivial approximately convex function defined on ``some part'' of it.

It remains to establish the connection with the normed structure
of the space. To this end, assume that $X$ is an ordered Banach space
\cite{ltz} and let $X^+=\{x\in X: x\geq 0\}$ be its positive cone. 
The continuous analogue of the Cholewa-Kominek function is
$$ x\longmapsto -\log_2\|x\|. $$

\begin{lemma}\label{-logx}
Let $X$ be an ordered Banach space. Then $-\log_2\|\cdot\|$ is $1$-convex
on $X^+\backslash \{0\}$.
\end{lemma}

\begin{proof} If $x$ and $y$ are positive,
then $\|x\|\leq \|x+y\|$. Thus, for $0\leq t\leq 1$, one has $\|t
x\| \leq \|t x + (1-t)y\|$ and so $\log_2 t + \log_2\|x\| \leq
\log_2\|t x + (1-t)y\|$; also $\log_2(1-t) + \log_2 \|y\| \leq
\log_2 \|t x + (1-t) y\|$. Therefore
$$
t \log_2 t +
(1-t) \log_2 (1- t) +t \log_2 \|x\| + (1-t)\log_2 \|y\|
\leq \log_2 \| t x + (1-t)y\|.
$$
Since $t \log_2t +
(1-t) \log_2 (1- t) \geq -1$ for all $0\leq t\leq1$ (the minimum value is
attained at $t=1/2$) the function 
$-\log_2\|\cdot \|$ is 1-convex on $X^+\backslash\{0\}$
\end{proof}

\begin{wrapfigure}[17]{r}{7.7cm}
\includegraphics{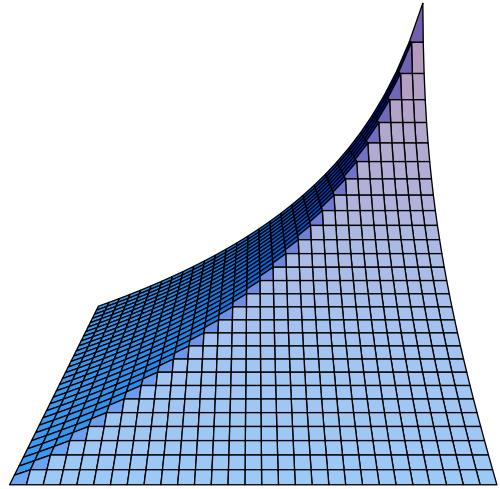}
\caption{
{\footnotesize 
The graph of the function $-\log_2\|\cdot\|$ on the positive part of the
ball of $\ell_\infty^2.$}
}
\end{wrapfigure}

\begin{lemma}\label{not}
With the above notations, if $X^+$ contains a weakly null
normalized sequence then $-\log_2\|\cdot\|$ cannot be uniformly
approximated by any convex function on $B_X^+\backslash\{0\}$.
\end{lemma}

\begin{proof}
Let $(u_n)$ be a weakly null sequence in $X^+$ with $\|u_n\|=1$.
Then there exists a sequence $(\sigma_n)$ of convex combinations of the
$u_n$'s such that $\|\sigma_n\|\to0$ and so $-\log_2\|(\sigma_n)\|\to \infty$.
However, any convex function approximating
$-\log_2\|\cdot\|$ should be bounded on convex combinations of the $u_n$ since
$-\log_2\|u_n\|=0$ for all $n$.
\end{proof}

Thus, we have defined a ``bad'' approximately convex function on
the unit ball of non-Schur cones. Observe
that in $\ell_1$ the preceding function becomes approximately convex.
Note, however, that composing
with the inclusion of $\ell_1$ into $c_0$ one obtains another
explicit ``bad'' 1-convex function on the infinite dimensional
simplex:
\begin{center}
$
(t_n)\in\Delta^\infty\longmapsto -\log_2\left(\max_n t_n\right).
$
\end{center}
To obtain the function defined on the whole unit ball requires
some extra work using additional properties of the spaces.

\begin{proposition} There exist nontrivial approximately convex
functions defined on the (closed) unit ball of either
$\ell_\infty$ or $ c$ (the space of convergent sequences).
\end{proposition}

\begin{proof}
Let $X$ denote one of the spaces $\ell_\infty$ or $c$. The key
point is that the unit ball of $X$ is affine-isomorphic to its
positive part: consider the map $B_{X}\to B^+_{X}$ given by
$x\longmapsto (1+x)/2$. Thus, it suffices to get a ``bad''
approximately convex function on the positive part of the unit
ball.

First of all, note that Lemma \ref{not} asserts that $x\longmapsto
-\log_2\|x\|$ is a ``bad'' $1$-convex map on
$B^+_{X}\backslash\{0\}$.

It remains to avoid the singularity at the origin.
Consider the ``involution'' of $B^+_{X}$ given by $x\longmapsto
1-x$ and let us modify the Cholewa-Kominek function to make $1$
the singular point instead of $0$. This yields the amended
function $$F(x) = - \log_2\|1 -x\|.$$ Now, for each $0<\theta<1$,
we set the 1-convex functions $$F_\theta (x) = - \log_2 \|1
-\theta x\|$$ defined on the whole $B^+_X$; they are increasingly
far from convex functions as $\theta$ approaches $1$. It remains
to past these pieces together.

Let $P_n$ be the projection of $X$ onto its first
$n$-coor\-dinates, and let $(\theta_n)$ be a sequence converging
to 1. We define the functions
$$F_n (x) = - \log_2\|P_n(1) -\theta_n P_n(x)\|$$
which are
non-negative and 1-convex on the positive part of the unit ball of
$X$.

It is clear that if $x_j=0$ for all $j\geq m$, then $F_j(x)=0$
for all $j\geq m$. It follows that for every finitely supported
$x\in B^+_{X}$ the sequence $F_n(x)$ is bounded.

Hence we can define a new $1$-convex function on the finitely
supported sequences of $B^+_{X}$ by taking $$ F^*(x)=\sup_n
F_n(x). $$ (The pointwise supremum of $\e$-convex functions is
$\e$-con\-vex whenever it makes sense.) Finally, we extend $F^*$
to all of $B^+_{X}$ by putting $F^*(x)=0$ for $x\notin c_{00}$.
Using the fact that the complement of $c_{00}$ in $B^+_{X}$ acts
as an ideal (with respect to convex combinations)  and that $F^*$
is non negative on $c_{00}\cap B^+_{X}$ it is easily verified that
the resulting map $F^*$ is $1$-convex on the whole of $B^+_{X}$.

The construction concludes showing that $F^*$ is at infinite
distance from any convex function on $B^+_{X}$. Actually, the
ensuing argument shows that $$ d_{c_{00}\cap
B^+_{X}}(F^*,g)=\infty $$ for all convex $g:c_{00}\cap
B^+_{X}\to\R$.

To see this, fix $n\geq 1$ and let $S_n=\sum_{i=1}^ne_i$. Consider
the points $p_i=S_n-e_i$. We see that $F^*(p_i)=0$. However,
$$
\frac{1}{n}\sum_{i=1}^np_i =\frac{n-1}{n}\cdot S_n 
$$
and so
$$\limsup_{n \to \infty} F^* \left( \frac{1}{n} \sum_{i=1}^n p_i
\right) \geq \limsup_{n \to \infty}
 \left(- \log_2 \left( 1 - \theta_n \frac{n -1}{n} \right)\right) =
+\infty$$
which already implies that $F^*$ cannot be approximated by any
convex function defined on a convex set containing $\{S_n-e_i:
n\in\mathbb N, 1\leq i\leq n\}$.
\end{proof}

We have thus solved the main problem in \cite{cabe} showing that
non-trivial $\e$-convex functions can occur in the unit ball of
$K$-spaces.

As for the other type of $K$-spaces, the B-convex ones, it
remains to know:

\begin{question}\label{bconvex}
\emph{Do B-convex spaces admit bad $\e$-convex functions on
their unit balls?}
\end{question}

\section{The stability of affinity in $K$-spaces}\label{sec:af}

We already know that the stability of convexity does not hold on
(the unit ball of) arbitrary $K$-spaces. However, under some
homogeneity conditions, $\e$-convex functions defined on
$K$-spaces can be approximated by convex ones (see, e.g.,
\cite[theorems 2 and 4]{cabe}). Thus one may wonder which other
classes of functions could provide a characterization of
$K$-spaces. Having in mind the behaviour of quasi-linear functions
on balls, it is not too surprising that $\e$-affinity works.
Recall that a mapping $f:D\to\R$ (here $D$ is a convex subset of a
linear space) is said to be $\e$-affine if it satisfies
\begin{equation*}
|f(tx+(1-t)y)-tf(x)-(1-t)f(y)|\leq\e
\end{equation*}
for all $x,y\in D,t\in[0,1]$. That is, an $\e$-affine map is one
which, in addition to being $\e$-convex, is also ``$\e$-concave''.
Before going further, let us prove the following simple consequence of
Hahn-Banach separation theorem.

\begin{proposition}
Let $f$ be any function defined on a convex set. 
If $f$ is approximable by a convex function and by a concave function, then
it is approximable by an affine function.
\end{proposition}

\begin{proof}
Put $\alpha=d_D(f,g), \beta=d_D(f,g)$, where $D$ is the underlying convex
set, $g$ is convex and $h$ concave. Then $-\alpha\leq f(x)-g(x)\leq\alpha$
and $-\beta\leq f(x)-g(x)\leq \beta$, hence
$$
-\beta+h(x)\leq f(x)\leq \alpha+g(x)\quad\quad(x\in D).
$$
But $-\beta+h$ is concave and $\alpha+g$ convex and so, if $X$ is any linear
space containing $D$, the sets
$$
E=\{(x,t): x\in D, t>\alpha+g(x)\}\quad\text{and}
\quad F=\{(x,t): X\in D, t<h(x)-\beta\}
$$
are nonoverlapping convex subsets of $X\times \R$.
The Hahn-Banach theorem gives an hyperplane $H$ separating 
$E$ from $F$, that is, an $\R$-linear functional
$a:X\times\R\to\R$ and $\gamma\in\R$ such that
$$
E\subset \{(x,t): a(x,t)\geq\gamma\}\quad\text{and}\quad
F\subset\{(x,t):a(x,t)\leq\gamma\}.
$$
Write $a(x,t)=b(x)+\lambda t$, where $b:X\to\R$ is $\R$-linear. It is clear
that $\lambda\neq 0$. It follows that $H$ is the graph of certain affine
function $a':X\to\R$ separating $-\beta+h$ from $\alpha+g$:
$$
-\beta+h(x)\leq a'(x)\leq \alpha+g(x)\quad\quad(x\in D).
$$
Hence $-2\beta\leq a'(x)-f(x)\leq 2\alpha$ and so $d_D(f,a')\leq
2\max\{\alpha,\beta\}$.
\end{proof}

In this way we obtain an alternate proof of the fact that
approximately affine functions are uniformly approximable by affine
functions on finite dimensional convex sets \cite[theorem~2]{lacz}. 
On the other hand, it is clear
that the restriction of a quasi-linear map to a bounded set is approximately
affine, and so is Ribe's original function on $\Delta^\infty$. Thus,
we obtain that the  ``affine'' version of Proposition~\ref{infinite} is also
true: if $X$ is an infinite dimensional Banach space, there is an
approximately affine function defined on a compact subset of $X$
which cannot be approximated by any affine function.

Our main result on $\e$-affine functions is the following theorem
about the local stability of affinity for $K$-spaces.

\begin{theorem}\label{staafin}
A Banach space $X$ is a $K$-space if and only if every $\e$-affine
function $f:B_X\to\R$ is (uniformly) approximable by a true affine
function $a:B_X\to \R$.
\end{theorem}

\begin{proof}
The `if' part is contained in Lemma \ref{midconvex}, taking into
account that if $f:X\to \R$ is quasi-linear with constant $\e$,
the restriction on the unit ball is $\e$-affine.

The `only if' part follows from the following slightly stronger
result which gives us more flexibility on the choice of the domain
$D$.
\end{proof}

\begin{proposition}\label{flex}
Let $D$ be a convex bounded set with nonempty interior in a
$K$-space $X$. Then to every $\e$-affine function $f:D\to\R$ there
corresponds an affine function $a:D\to\R$ such that $$
d_D(f,a)\leq A\cdot \e, $$ where $A=A_D$ is a constant depending
only on $D$.
\end{proposition}

\begin{proof}
There is no loss of generality in assuming that the origin is
interior to $D$ and also that $f$ is 1-affine, with $f(0)=0$. 
Let $\mathscr L$ denote the set
of lines of $X$ passing throughout $0$ and, for each
$\ell\in\mathscr L$, set $D_\ell=D\cap\ell$ and let $f_\ell$ be
the restriction of $f$ to $D_\ell$.

Each $f_\ell$ is clearly $\e$-affine and since $D_\ell$ is
one-dimen\-sional, there is an affine function
$a_\ell:D_\ell\to\R$ with $$ d_{D_\ell}(f_\ell,a_\ell)\leq 1. $$
(See \cite{lacz}.) In particular, we have $|a_\ell(0)|\leq 1$ and
so we can extend $x\mapsto a_\ell(x)-a_\ell(0)$ to a linear map
$L_\ell:\ell\to\R$ with $$ d_{D_\ell}(f_\ell,L_\ell)\leq 2. $$ Let
us define $f^*:X\to\R$ by $$ f^*(x)=L_{[x]}(x), $$ where $[x]$ is
the line spanned by $x$. It is clear that $f^*$ is homogeneous on
$X$ and also that $d_D(f,f^*)\leq 2$, from where it follows that
$f^*$ is $3$-affine on $D$.

We see that $f^*$ is quasi-linear. Indeed, for $x,y\in D$ one has
$$
|f^*(x+y)-f^*(x)-f^*(y)|
=2\left|f^*\left(\frac{x+y}{2}\right)-\frac{f^*(x)+f^*(y)}{2}\right|\leq 6.
$$
It follows that $f^*$ is quasi-additive, with
$Q(f^*)\leq 6/r_0$, where $r_0=\sup\{r>0:rB_X\subset D\}$.
Now, since $X$ is a $K$-space, there is a linear (hence affine)
function $a:X\to R$ satisfying $$ |f^*(x)-a(x)|\leq
\frac{6M}{r_0}\|x\|\quad\quad(x\in X), $$ where $M$ depends only
on $X$. Therefore, if $R_0=\inf\{R>0:D\subset RB_X\}$, we have
$d_D(f^*,a)\leq 6MR_0/r_0$ and so
\begin{equation*}
d_D(f,a)\leq \frac{6MR_0}{r_0}+2,
\end{equation*}
which completes the proof.
\end{proof}

\section{Local stability of the Jensen equation}

Suppose $D$ is a midpoint convex subset of a linear space.
Jensen's functional equation is $$
f\left(\frac{x+y}{2}\right)=\frac{f(x)+f(y)}{2}\quad\quad(x,y\in
D). $$ Accordingly, we say that $f$ is $\e$-Jensen if $$
\left|f\left(\frac{x+y}{2}\right)-\frac{f(x)+f(y)}{2}\right|\leq\e
\quad\quad(x,y\in D). $$ Let us prove the Jensen analogue of
Theorem \ref{staafin}:

\begin{theorem}\label{stajens}
A Banach space $X$ is a $K$-space if and only if every $\e$-Jensen
function $f:B_X\to\R$ is (uniformly) approximable by a function
$a:B_X\to \R$ which satisfies Jensen equation.
\end{theorem}

The `if' part follows from Lemma \ref{midconvex}, while the
converse is contained in the following result.

\begin{proposition}\label{flex-j}
Let $D$ be a convex bounded set with nonempty interior in a
$K$-space $X$. Then for every $\e$-Jensen function $f:D\to\R$
there is a true Jensen function $a:D\to\R$ such that $
d_D(f,a)<\infty$.
\end{proposition}

\begin{proof} The proof is almost the same as before and we only
give the main steps. Assume again that $0$ is interior to $D$ and
that $f:D\to\R$ is 1-Jensen, with $f(0)=0$. With the same notations, 
it is clear that every $f_\ell$
is 1-Jensen and since $D_\ell$ is one-dimensional and convex
there is a Jensen function $a_\ell:D_\ell\to\R$ with $$
d_{D_\ell}(f_\ell,a_\ell)\leq 2. $$ (See \cite[theorem 3]{lacz}.)
Putting $L_\ell(x)=a_\ell(x)-a_\ell(0)$ we obtain a 2-homogeneous
Jensen function that clearly extends to a 2-homogeneous Jensen
function on $\ell$ (which do not relabel), with $
d_{D_\ell}(f_\ell,L_\ell)\leq 4$. Setting $$ f^*(x)=L_{[x]}(x), $$
(where $[x]$ is the line spanned by $x$), it is clear that $f^*$
is $2$-homogeneous on $X$ and also that $d_D(f,f^*)\leq 4$, from
where it follows that $f^*$ is $6$-Jensen on $D$.

Thus, for $x,y\in D$, one has
$$
|f^*(x+y)-f^*(x)-f^*(y)|
=2\left|f^*\left(\frac{x+y}{2}\right)-\frac{f^*(x)+f^*(y)}{2}\right|\leq 12.
$$
 Hence $f^*$ is quasi-additive, with
$Q(f^*)\leq 12/r_0$, where $ r_0=\sup\{r>0:rB_X\subset D\} $.

Now, since $X$ is a $K$-space, we can use \cite[corollary 2]{sing}
to get an additive (hence Jensen) function $a:X\to \R$ with
$f^*-a$ continuous at the origin of $X$. Therefore, there is
$\delta>0$ so that $$ |f^*(x)-a(x)|\leq
1\quad\quad(\|x\|\leq\delta). $$ Having in mind that both $f^*$
and $a$ are 2-homogeneous we see that, in fact, $$
|f^*(x)-a(x)|\leq\frac{\|x\|}{\delta}\quad\quad(x\in X). $$ In
particular one has $d_D(f^*,a)\leq R_0/\delta$, where
$R_0=\inf\{R>0:D\subset RB_X\}$ and so
\begin{equation}\label{so}
d_D(f,a)\leq
4+\frac{R_0}{\delta},
\end{equation}
which completes the proof.
\end{proof}

Following Laczkovich \cite{lacz}, given a convex set $D$, we write
$C_D$ (respectively, $A_D$ and $J_D$) for the least constant
$\kappa$ such that to every $\e$-convex (respectively, $\e$-affine
and $\e$-Jensen) function $f:D\to R$ there corresponds a convex
(respectively, affine and Jensen) function $g:D\to \R$ such that
$d_D(f,g)\leq \kappa\e$.

One of the most surprising results in \cite{lacz} is that (in the
finite dimensional case) $C_D$ is essentially independent on the
shape of $D$. In fact $C_D\sim \log\dim D$, where $\dim D$ is the
linear dimension of the least affine submanifold containing $D$.

Laczkovich also observed that both $A_D$ and $J_D$ are of the same
order as $\log\dim D$ if $D$ is either a simplex or an
$n$-dimensional ``octahedron'' (the unit ball of $\ell_1^n$). He
 asks for the behaviour of the constants $A_D$ and
$J_D$ for $D$ either the $n$-dimensional cube (the unit ball of
$\ell_\infty^n$) or the $n$-dimensional euclidean ball (the unit
ball of $\ell_2^n$). It turns out that these constants are
uniformly bounded:

\begin{corollary}
If $D$ is the unit ball of one of the spaces $\ell_\infty^n,c_0$
or $ \ell_\infty$ (or even of an $\mathcal L_{\infty,1}$-space),
then $A_D\leq 6\cdot 200+2$ and $J_D\leq 2\cdot(6\cdot 200+2)$.

If $D$ is the unit ball of a Hilbert space, then $A_D\leq 6\cdot
37+2$ and $J_D\leq 2\cdot(6\cdot 37+2)$.
\end{corollary}

\begin{proof}
The statement concerning the constants $A_D$ in the first case
follows from the fact (proved by Kalton and Roberts
\cite{kaltrobe}) that every $\mathcal L_{\infty,1}$-space is a
$K$-space, with constant at most $200$ (although 100 was announced
in \cite{kaltrobe}; see \cite{MR}). The estimate for $J_D$ follows
from the fact that $J_D\leq 2A_D$ for finite dimensional $D$ (see
\cite{lacz}) and
local techniques: we give only a sketch.
Let $\mathscr F$ be the set of all finite
dimensional subspaces of $X$.
Since $X$ is an $\mathcal L_{\infty,1}$-space, given $F\in\mathscr
F$ and $\e>0$, there is $E\in\mathscr F$ containing $F$ and such
that the Banach-Mazur distance between $E$ and $\ell_\infty^{\dim
E}$ is at most $1+\e$. Which simply means that there is a linear
isomorphism $T:E\to\ell_\infty^{\dim E}$ such that $$
B_{\ell_\infty^{\dim E}}\subset TB_E\subset (1+\e)
B_{\ell_\infty^{\dim E}}. $$ Thus, in view of (\ref{so}), one has
$$ A_{B_E}\leq 6\cdot 200\cdot(1+\e)+2 $$ and therefore $$
J_{B_E}\leq 2\cdot(6\cdot 200\cdot(1+\e)+2) $$ Since $\e$ is
arbitrary, we see that if $f$ is $1$-Jensen on $B_X$ then for
every $F\in \mathscr F$ there is a Jensen $a_F:B_F\to\R$ such that
$$ d_{B_F}(f,a_F)<2\cdot(6\cdot 200+2)+\frac{1}{\dim F} $$ To end,
let $\mathfrak V$ be an ultrafilter refining the Fr\'echet
(=order) filter on $\mathscr F$ and define $a:B_X\to\R$ taking $$
a(x)=\lim_{\mathfrak V(F)}a_F(x). $$ It is easily seen that the
above definition yields a Jensen function at distance at most
$2\cdot(6\cdot200+2)$ from $f$ on $B_X$.

As for euclidean norms, the standard proofs that Hil\-bert spaces
are $K$-spaces do not give any estimate for the corresponding
constant \cite{kaltf,kspace}. There is however a recent paper by
\v Semrl \cite{semrl} in which it is shown that if $f:H\to\R$ is
quasi-linear and bounded on the unit ball, then there exists a
linear map $\ell:H\to\R$ such that $$ \dist(f,\ell)\leq 37\cdot
Q(f). $$ Since quasi-linear maps are always bounded on finite
dimensional balls an obvious local argument shows that Hil\-bert
spaces are $K$-spaces with constant 37. The result follows from
this.
\end{proof}

It will be clear for those acquainted with twisted sums of Banach
spaces that the results of sections 3 and 4 cannot be extended to
vector valued maps. If $X$ is any infinite dimensional Banach
space there is another Banach space $Y$ and a quasi-linear map
$f:X\to Y$ such that $\dist(f,\ell)=\infty$ for all linear maps
$\ell:X\to Y$.

(Indeed, if $X$ contains a complemented subspace isomorphic to
$\ell_1$, set $Y=\R$ and use Ribe's map in the obvious way.
Otherwise, take a quotient operator $\pi:\ell_1(\Gamma)\to X$ for
a suitable set $\Gamma$, set $Y=\ker\pi$ and note that the exact
sequence
$$
0\to Y\to\ell_1(\Gamma)\to X\to 0
$$
does nor split. This implies the existence of a quasi-linear map
$f:X\to Y$ at infinite distance from all linear maps $X\to Y$; see
\cite{kaltf} or \cite{diss}.)

It is clear that $f:B_X\to Y$ is approximately affine, but no
Jensen map $a:B_X\to Y$ is uniformly close to $f$. For if $a$ is a
Jensen map such that $d_{B_X}(f,a)<\infty$ then $a$ is actually
affine and so $f$ is asymptotically close to the linear map
$\ell:X\to Y$ obtained by extending (by homogeneity) $x\to
a(x)-a(0)$ to all of $X$.

\section{Uniform boundedness}\label{sec:unif}

In this Section we put the results of the preceding two in its proper
setting by showing that, for many convex sets $D$, if every
$\e$-affine (respectively, $\e$-Jensen) function $f:D\to\R$ is approximable
by an affine (respectively, Jensen) $a:D\to\R$, then this can be achieved
with $d_D(f,a)\leq C\e$, where $C$ is a constant depending only on $D$.

So, let $D$ be a convex set, where no topology is assumed. A point $d\in D$
is said to be geometrically interior to $D$ if for every $v$ in the
linear space spanned by $D$ one has $d\pm tv\in D$ for $t>0$ small enough.
This is equivalent to the following: $D$ is a neighbourhood of $d$
in the linear space spanned by $D$ equipped with the strongest locally convex
topology.

In the sequel, we say that $D$ is thick if it has at least one
geometrically interior
point. Every set with nonempty interior in a normed space is thick, and so
are all separable polytopes (a polytope is a closed bounded set in a Banach
space whose finite dimensional sections have all finitely many extreme
points; see \cite[\S 6]{fonflind}).

As for elementary simplices, 
let $\mathfrak c$ an infinite cardinal (that we regard also as an index
set). Then the $\mathfrak c$-dimensional simplex
$$
\Delta^\mathfrak c=\left\{x\in\R^\mathfrak c: x_\alpha\geq 0 \text{ for all }
\alpha\in\mathfrak c \text{ and } \sum_{\alpha\in\mathfrak
c}x_\alpha=1\right\}.
$$
is thick if and only if $\mathfrak c$ is countable.

\begin{theorem}\label{th:unif-j}
Let $D$ be a thick convex set. If every approximately Jensen function
on $D$ is uniformly approximable by some Jensen function, then there is a
constant $J$, depending only on $D$, such that, to every $\e$-Jensen
$f:D\to\R$ there corresponds a Jensen function $a:D\to\R$ with $d_D(f,a)\leq J\cdot
\e$.
\end{theorem}

This Theorem applies to the sets appearing in Proposition~\ref{flex-j}, and
thus $J_D$ is finite (as we already know is $A_D$) 
if $D$ is a bounded convex set with nonempty interior in a K-space.

The main step in the proof of Theorem~\ref{th:unif-j}
is the next Lemma showing that if $f$ is approximately Jensen, then a
judicious choice of the affine function $a:D\to\R$ allows us to control
$|f(x)-a(x)|$ ``pointwise'' on $x$ but ``uniformly'' on $f$.

\begin{lemma}\label{judicious}
Let $D$ be a thick set. There is a function $\eta:D\to\R$ (depending only on
$D$) such that, for every $\e$-Jensen function $f:D\to\R$ there is a Jensen
$a:D\to\R$ satisfying the estimate
$$
|f(x)-a(x)|\leq \e\eta(x)\quad\quad(x\in D).
$$
\end{lemma}

\begin{proof}
Without loss of generality we assume that the origin is geometrically
interior to $D$. Let $X$ be the linear space spanned by $D$ and consider the
gauge of $D$, that is, the function $\vr:X\to[0,\infty]$ defined as
$$
\vr(x)=\inf\{\lambda>0:\lambda^{-1}x\in D\}.
$$
The hypothesis on $D$ implies that
$\vr$ is well-defined and also that it takes only finite values.

Now, let $f$ be 1-Jensen on $D$ with $f(0)=0$. Just as in the first part of
the proof of Proposition~\ref{flex-j} we can find a function $f^*:X\to\R$
such that:
\begin{itemize}
\item $d_D(f,f^*)\leq 4$, and
\item $f^*$ is homogeneous over the rationals: $f^*(qx)=qf^*(x)$ for all
$x\in X$ and $q\in\mathbb Q$.
\end{itemize}
It follows that $f^*$ is 6-Jensen on $D$ and so
$$
|f^*(x+y)-f^*(x)-f^*(y)|\leq 12\quad\quad(x,y\in D).
$$
Using the $\mathbb Q$-homogeneity of $f^*$ we obtain that $f^*$ is
quasi-additive with respect to $\vr$:
$$
|f^*(x+y)-f^*(x)-f^*(y)|\leq 12(\vr(x)+\vr(y))\quad\quad(x,y\in X).
$$
It follows (by induction on $n$; see \cite{kaltf}) that
\begin{equation}\label{kal}
\left|f^*\left(\sum_{i=1}^nx_i\right)-\sum_{i=1}^nf^*(x_i)\right|\leq
12\sum_{i=1}^ni\vr(x_i)\quad\quad(n\in\N,x_i\in X).
\end{equation}
Finally, let $\mathscr B$ be a Hamel basis of $X$ over the rationals, with
$\pm b\in D$ for all $b\in\mathscr B$ and define a $\mathbb Q$-linear 
(hence Jensen) map $a:X\to\R$ taking $a(b)=f^*(b)$ for
$b\in\mathscr B$
and extending linearly on the rest.

Fixing $x\in X$, let us estimate $f^*(x)-a(x)$. Write $x=\sum_{b\in\mathscr B}
q_b b$ and let $n(x)$ be the number of nonzero summands in that
decomposition, so that we can write
$x=\sum_{i=1}^{n(x)}q_{b(i)}b(i)$. We have
\begin{align*}
|f^*(x)-a(x)|&=\left|f^*(x)-
\sum_{b\in\mathscr B} q_b f^*(b)\right|
=\left|f^*\left(\sum_{i=1}^{n(x)} q_{b(i)} b(i)\right)-
\sum_{i=1}^{n(x)} f^*(q_{b(i)}b(i)\right|\\
&\leq 12\sum_{i=1}^{n(x)}i\vr(q_{b(i)}b(i))\leq
12n(x)\sum_{b\in\mathscr B}|q_b|.
\end{align*}
Therefore, for $x\in D$ we get
$$
|f(x)-a(x)|\leq 4+12n(x)\sum_b|q_b|,
$$
and choosing $\eta(x)=4+12n(x)\sum_b|q_b|$ we conclude the proof.
\end{proof}

\begin{proof}[End the proof of Theorem~\ref{th:unif-j}] Let
$\mathscr{AJ}(D)$ denote the linear space of all approximately Jensen functions on $D$
and $\mathscr J(D)$ that of Jensen functions. We introduce two functionals 
on
$\mathscr{AJ}(D)$ as follows:
\begin{align*}
\e_J(f)&=\sup_{x,y\in
D}\left|f\left(\frac{x+y}{2}\right)-\frac{f(x)+f(y)}{2}\right|=\inf\{\e: f
\text{ is $\e$-Jensen}\}\\
\delta_J(f)&=d_D(f,\mathscr J(D)).
\end{align*}

With these notations our hypothesis is nothing but $\delta_J(f)<\infty$ for all
$f\in\mathscr{AJ}(D)$ and we must prove that $\delta_J(f)\leq J\cdot\e_J(f)$
for some constant $J$ independent on $f$. Obviously, $\e_J(\cdot)\leq
2\delta_J(\cdot)$.

It is clear that both $\e_J$ and $\delta_J$ are
seminorms. We have
$\ker \e_J=\ker \delta_J=\mathscr J(D)$, so that both
$\e_J$ and $\delta_J$ are well-defined norms on the quotient
space $\mathscr{AJ}(D)/\mathscr J(D)$.
That $\mathscr{AJ}(D)/\mathscr J(D)$ is complete under $\delta_J$ is nearly
obvious. Thus, the following result finishes the proof, thanks to the open
mapping theorem.\end{proof}

\begin{lemma}
The space $\mathscr{AJ}(D)/\mathscr J(D)$ is complete under the norm $\e_J$.
\end{lemma}

\begin{proof}
It suffices to show that absolutely summable series
converge in $\mathscr{AJ}(D)/\mathscr J(D)$.
So let $(f_n)$ be such that
$$
\sum_{n=1}^\infty \e_J(f_n)<\infty.
$$
By Lemma~\ref{judicious} there are Jensen functions $a_n:D\to\R$ such that
$
|f_n(x)-a_n(x)|\leq\e_J(f_n)\eta(x)
$
 for all $x\in D$. Hence we can define a function $f$ pointwise as
$$
f(x)=\sum_{n=1}^\infty (f_n(x)-a_n(x)).
$$
It is straightforward that $g$ is approximately Jensen, with $\e_J(f)\leq
\sum_n\e_J(f_n)$ and also that $\sum_n[f_n]$ converges to $[f]$ in
$\mathscr{AJ}(D)/\mathscr J(D)$ with the norm $\e_J$.
\end{proof}

The following ``affine'' companion of Theorem~\ref{th:unif-j} has a simpler
proof we leave to the reader.

\begin{theorem}\label{th:unif-a}
Let $D$ be a thick convex set. If every approximately affine function
on $D$ is uniformly approximable by some affine function, then there is a
constant $A_D$, such that, to every $\e$-affine
$f:D\to\R$ there corresponds an affine $a:D\to\R$ with $d_D(f,a)\leq A_D\cdot
\e$.\hfill $\square$
\end{theorem}

We do not know if the hypothesis on $D$ can be removed from the results in
this Section. Also, it would be interesting to get similar results for
approximately convex functions. However we strongly believe that if every
approximately convex function on $D$ is approximable, then $D$ has finite
dimension. If so, the corresponding uniform boundedness result for
approximate convexity would be a trivial tautology.

\section{$c_0$ is not isometric to a Banach envelope \ldots}

Let us recall from \cite{kaltpeckrobe} the minimal background
that one needs to understand what follows. A quasi-norm on a
(real or complex) vector space $X$ is a non-negative real-valued
function on $X$ satisfying:
\begin{itemize}
\item  $\|x\|=0$ if and only if $x = 0$;
\item  $\|\lambda x\|=|\lambda|\|x\|$ for all $x\in X$ and
$\lambda\in\mathbb K$;
\item  $\|x+y\|\leq \Delta(\|x\|+\|y\|)$ for some fixed $\Delta\geq
1$ and all $x, y \in X$.
\end{itemize}
A quasi-normed space is a vector space $X$ together with a
specified quasi-norm. On such a space one has a (linear) topology
defined as the smallest linear topology for which the set $B_X =
\{x \in X : \|x\|\leq 1\}$ (the unit ball of $X$) is a
neighborhood of 0. In this way, $X$ becomes a locally bounded
space (that is, it has a bounded neighborhood of 0); and,
conversely, every locally bounded topology on a vector space comes
from a quasi-norm. A quasi-Banach space is a complete quasi-normed
space.

Needless to say, every Banach space is a quasi-Banach space, but
there are important examples of quasi-Banach spaces which are not
(isomorphic to) Banach spaces. Let us mention the $\ell_p$ and
$L_p$ spaces and the Hardy classes $H^p$ for $0<p<1$.

Let $X$ be a quasi-Banach space. The dual space $X^*$ is always a
Banach space under the norm
$$
\|x^*\|=\sup_{\|x\|\leq1}|x^*(x)|.
$$
Consider the ``evaluation mapping" $\delta : X\to X^{**}$ given by
$\delta(x)x^* = x^*(x)$. The Banach envelope $\co X$ of $X$ is the
closure of $\delta(X)$ in $X^{**}$ equipped with the induced norm.
Notice that $\delta$ is one to one if and only if $X^*$ separates
$X$ in the sense that for every nonzero $x\in X$ there is $x^*$ in
$X^*$ such that $x^*(x)\neq 0$.  So, $\co X$ is a Banach space
whose unit ball equals the closed convex hull of $\delta(B_X)$. In
particular, when $X$ is finite-dimensional $\co X$ can be seen as
a renorming of $X$ itself. The map $\delta : X \to \co X$ has the
following universal property: every bounded linear operator from
$X$ into a Banach space $Y$ factorizes throughout $\delta$ with
equal norm. This clearly implies that $\co X$ is the ``nearest"
Banach space to $X$ with respect to the Banach-Mazur distance: if
$T$ is an isomorphism from $X$ into a Banach space $Y$, then
$\delta:X\to \co X$ is an isomorphism and, in fact
$\|T\|\|T^{-1}\|\leq\|\delta\|\|\delta^{-1}\|$. Note that
$\|\delta\|=1$ and that $\|\delta^{-1}\|$ is the least constant
$K$ for which
$$
\|\delta  x\|_{\co X}\leq\|x\|_X\leq K\|\delta x\|_{\co X}
$$
holds for all $x\in X$. Of course, $X$ is locally convex (that is,
isomorphic to a Banach space) if and only if $\|\cdot\|_{\co X}$
is equivalent to the original norm $\|\cdot\|_X$.

From the point of view we have adopted  in this paper, it is worth
noting that $\|\cdot\|_{\co X}$ equals $\co\|\cdot\|_X$ (on the
common domain $X$) and that the Lipschitz counterexamples found in
\cite{casipapi}) depend on the fact that $\ell_1$ is the Banach
envelope of the non-locally convex spaces $\ell_p$ for $0<p<1$.

Which Banach spaces can be envelopes of non-locally convex spaces
with separating dual? (The condition of having separating dual is
to avoid trivial examples.) In \cite{kaltf}, Kalton proved that
the Banach envelope of a non-locally convex quasi-Banach space
with separating dual is never B-convex. He then asks whether $c_0$
(or $l_\infty$) can be the Banach envelope of a non-locally convex
quasi-Banach space with separating dual \cite{kaltf}. Kalton
himself solved the ``isomorphic'' problem in the negative
\cite[section 4]{combo} showing a rather pathological quasi-Banach
space whose Banach envelope is isomorphic to $c_0$. Surprisingly
one has:

\begin{proposition}
\label{isom}
$c_0$ is not isometric to the
Banach envelope of a nonlocally
convex quasi-Banach space with separating dual.
\end{proposition}

This straightforwardly follows after the following result.

\begin{theorem}\label{c0}
Let $X$ be a quasi-Banach space. If $T:X\to c_0$ is a bounded operator such
that $\overline{\co} TB_X=B_{c_0}$ then $T$ is an open mapping of $X$ onto
$c_0$.
\end{theorem}

The proof is based of a few elementary observations that we put
together in the following Lemma. Let us say that $u\in B_{c_0}$ is
a locally extreme point if for all $k$ one has $|u_k|\in\{0,1\}$ .

\begin{lemma}\label{ju}
With the same notations as in Theorem \ref{c0}.
\begin{itemize}
\item[(a)] For every $y\in B_{c_{00}}$ there is a locally
extreme point $u\in B_{c_0}$ such that $\|y-\frac{1}{2}u\|\leq \frac{1}{2}$.
\item[(b)] For every locally extreme
point $u$ and every $\varepsilon > 0$, there is $x\in B_X$ such
that $\|u-1_{\supp(u)}T(x)\|\leq\varepsilon$.
\item[(c)] For every $y$ in the unit ball of $c_0$ and every $\e>0$ there is
$x\in B_X$ such that $\|y-\frac{1}{2}T(x)\|\leq \frac{1}{2}+\e$.
\end{itemize}
\end{lemma}

\begin{proof}
(a) Take $u=\sgn y$.

(b) Let $u$ be a locally extreme point of the ball of $c_0$ and
consider the (contractive) projection $1_{\text{supp}(u)} : c_0\to
c_0$ given by multiplication. The hypothesis implies that
$\overline \co (1_{\text{supp}(u)}TB_X) =
B(1_{\text{supp}(u)}c_0)$. Since $1_{\text{supp}(u)}c_0$ is finite
dimensional, the closure of $1_{\text{supp}(u)}TB_X$ contains all
extreme points of the unit ball of $1_{\supp(u) }c_0$. In
particular, it contains $u$.

(c) We may assume $y\in c_{00}$. Let $u=\sgn y$ and note that
$\supp(u)=\supp(y)$. By (b)  there is $x\in B_X$ such that
$\|u-1_{\supp(u)}T(x)\|\leq\varepsilon$. One thus has
$$
\Bigl{|}y(k)-\frac{1}{2}Tx(k)\Bigr{|}\leq
\begin{cases}
\frac{1}{2}+\varepsilon
&k\in\supp(f)\\
\frac{1}{2}
 &k\notin\supp(f),
 \end{cases}
$$
and the result follows.
\end{proof}

\begin{proof}[Proof of Theorem \ref{c0}]
Fix $y\in c_0$, with $\|y\|\leq 1$ and let $\e>0$ be fixed. By
Lemma \ref{ju} there is  $x_1\in B_X$ such that
$\|y-\frac{1}{2}Tx_1\|\leq \frac{1}{2}+\e$. Replacing $y$ by
$y-\frac{1}{2}Tx_1$, there is $x_2\in B_X$ such that
$$
\left\|y-\frac{1}{2}Tx_1-\frac{1}{2}\left(\frac{1}{2}+\e\right)Tx_2\right\|\leq
\left(\frac{1}{2}+\e\right)^2.$$ Proceeding inductively one
obtains a sequence $(x_n)$ in the unit ball of $X$ such that
\begin{equation}\label{that}
 \left\|y-\frac{1}{2}\sum_{i=1}^n\left(\frac{1}{2}+\e\right)^{i-1}Tx_i\right\|\leq
\left(\frac{1}{2}+\e\right)^n
\end{equation}
holds for all $n$. Let us estimate the quasi-norm of
$\sum_{i=1}^n\left(\frac{1}{2}+\e\right)^{i-1}x_i.$

By the Aoki-Rolewicz
theorem \cite{kaltpeckrobe}, $X$ has an
equivalent $p$-norm 
for some $0<p\leq 1$. Hence, for
some constant $M=M(p,X)$ and all $n$ one has
$$
\left\|\sum_{i=1}^n\left(\frac{1}{2}+\e\right)^{i-1}x_i\right\|^p
\leq
M
\sum_{i=1}^n\left\|\left(\frac{1}{2}+\e\right)^{i-1}x_i\right\|^p
= M \sum_{i=1}^n\left(\frac{1}{2}+\e\right)^{(i-1)p}
\leq \frac{M}{1-\left(\frac{1}{2}+\e\right)^p}
$$
for $\e$ small enough. Thus
$$
\left\|\sum_{i=1}^n\left(\frac{1}{2}+\e\right)^{i-1}x_i\right\|\leq
K
$$
for some constant $K$. Taking (\ref{that}) into account, it is
clear that $T$ is almost open. An appeal to the (proof of the)
open mapping theorem \cite{kaltpeckrobe} ends the proof.
\end{proof}

We can use the above argument to show that B-convex spaces cannot
be envelopes of non-locally convex spaces with separating dual.
First of all, note that every quasi-norm is equivalent to some
$p$-norm for some $0<p\leq 1$ (this is the Aoki-Rolewicz theorem).
This implies that if $\co Y=X$ isometrically and $Y$ has an
equivalent $p$-norm, then $X=\co \left(\co_p Y\right)$, where
$\co_p Y$ denotes the $p$-Banach envelope of $Y$ (here the
$p$-convex envelope of a symmetric set $A$ is defined 
to be the set $ \co_p
A= \left\{\sum_{i=1}^nt_ia_i : a_i\in A, \sum_{i=1}^n|t_i|^p\leq 1\right\}$). 
And it is so because
$$
\overline \co B_Y = \overline \co \left(\overline \co_p
B_Y\right),
$$
Let us say that $B$ has the property $(p,\theta,\kappa)$ if
whenever $A$ is a symmetric $p$-convex subset of $B$ such that
$B=\overline\co A$ one has
$$
B\subset \theta B+\kappa A.
$$
For instance, the content of Lemma \ref{ju} is that the unit ball
of $c_0$ has all properties $(p,\frac{1}{2}+\e,2)$ for $0<p<1$ and
$\e>0$. Another non-trivial example is provided by the unit ball
of any B-convex space:

\begin{lemma}
Let $X$ be a B-convex Banach space. Then, for every $0<p<1$ there
are $\theta<1$ and $\kappa>0$ such that  $B_X$ has the property
$(p,\theta,\kappa)$.
 \end{lemma}

\begin{proof} This follows from Bruck's \cite[theorem 1.1]{bruck}
that identifies B-convexity with the so-called convex approximation
property. Namely, $X$ is B-convex if and only if given a
bounded $B\subset X$ and $\e>0$ there exists $m\in\N$ such that
$\co B\subset \co^{[m]} B+\e B_X$,
where $$ \co^{[m]} B=\left\{\sum_{i=1}^mt_ib_i: \sum_{i=1}^mt_i=1,
t_i\geq 0, b_i\in B\right\}.
$$
Now, if $A\subset B_X$ is a $p$-convex set such that
$B_X=\overline \co A$, we have
$
B_X \subset \co^{[m]} A+\e B_X
$
for all $\e>0$ and some $m \in \N$. Since $A$ is $p$-convex
there is $M>0$ such that $\co^{[m]} A \subset MA$ and so $B_X$ has
$(p,\e,M)$, as desired.
\end{proof}
The proof of the following result is contained in that of Theorem
3. We make it explicit for the sake of clarity.

\begin{proposition}
If $B_X$ has the property $(p,\theta,\kappa)$ for some $\theta <1$
then $X$ cannot be the Banach envelope of a nonlocally convex
$p$-Banach space with separating dual.\hfill$\square$
\end{proposition}

We obtain an alternate proof of Kalton result:

\begin{corollary}
$B$-convex spaces cannot be envelopes of non-locally convex spaces
with separating dual.\hfill$\square$
\end{corollary}

It is not true that admitting an equivalent norm with the
properties $(p, \theta, \kappa)$ guarantees that the space is a
$K$-space: it is well-known that $H^1$ is not a $K$-space (it
contains $\ell_1$ as a direct factor) yet its usual norm has all
properties $(p,\frac{1}{2},1)$ since $H^1$ has the Radon-Nikod\'ym
property and every norm one $f\in H^1$ can be written as the
midpoint of two inner functions (see \cite{du, katz}).

\section{\ldots  although it
is $(1+\e)$-isomorphic!}

In spite of the previous negative result, we show now that there are
renormings of $c_0$, in fact small perturbations of the original norm, that 
are envelopes of nonlocally convex spaces with separating dual.

\begin{wrapfigure}[19]{r}{7.9cm}
\includegraphics{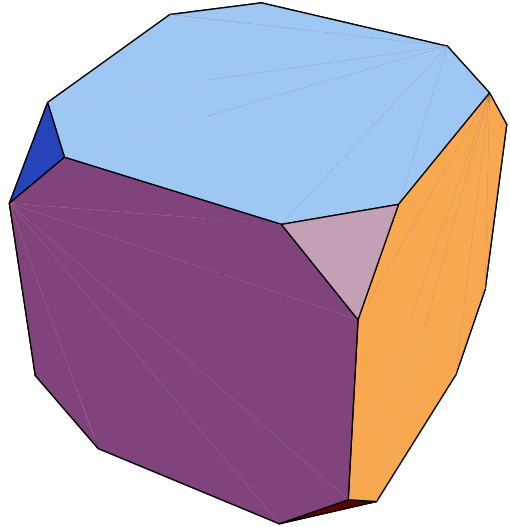}
\caption{
{\footnotesize 
The unit ball  of $\co X_p(\Omega)$ in three dimensions}
}
\end{wrapfigure}

This shows
that
the set of
Banach envelopes fails to be closed with respect to the
Banach-Mazur distance.

Precisely, we show:

\begin{proposition}
For every $\varepsilon>0$ there exists a non-locally convex
quasi-Banach space with separating dual whose Banach envelope
is $(1+\varepsilon)$-isomorphic to $c_0$.
\end{proposition}
\begin{proof} The construction is a  generalization of Kal\-ton's original
example (which is 2-isomorphic to $c_0$) and is based on (a refinement of)
an observation used by Talagrand \cite{tal} to construct a pathological
submeasure.

For each rational $\varepsilon>0$ and each positive integer $n$ for which
$n\varepsilon$ is an integer, consider the set
\smallskip

\begin{center}

$S(\varepsilon,n)=\{1,2,\dots,(1+\varepsilon)n\}$ 
\end{center}
\smallskip

\noindent
and let
$\Omega=\Omega(\varepsilon,n)$
be the class of 
all subsets of $S(\varepsilon,n)$ having
cardinal $n$.
Finally, set $A_i=\{\omega\in\Omega:i\in\omega\}$ for $1\leq i\leq
(1+\varepsilon)n$.

\begin{lemma}\label{talagrand}
The collection
$\{A_i\}_{i=1}^{(1+\varepsilon)n}$ has the following properties:
\begin{itemize}
\item[(a)] If $J\subset\{1,2,\dots,(1+\varepsilon)n\}$ is such that
$|J|\leq\varepsilon n$, then $\cup_{i\in J}A_i \neq \Omega$.
\item[(b)] $\sum_{i=1}^{(1+\varepsilon)n}1_{A_i}=n1_\Omega$.
\end{itemize}
\end{lemma}
\begin{proof}(a) Take an $\omega\in\Omega$ that does not intersect $J$ and
note that $\omega\in A_i$ if and only if $i\in\omega$. It is clear
that $\omega\notin\cup_{i\in J}A_i$.

(b) Let $\omega$ be arbitrarily chosen in $\Omega$. Then
$$
\sum_{i=1}^{(1+\varepsilon)n}1_{A_i}(\omega)
=\sum_{i=1}^{(1+\varepsilon)n}1_\omega(i)=|\omega|=n
$$
\end{proof}
Now fix $p$, with $0<p<1$, and define a quasi-norm on the space
$\ell_\infty(\Omega)$ by putting
$$
\|f\|_p=\inf\left\{\left(\sum_{i=1}^{(1+\varepsilon)n}|c_i|^p\right)^{1/p}:
|f|\leq\sum_{i=1}^{(1+\varepsilon)n}c_i1_{A_i}\right\}.
$$
Let $X_p(\Omega)$ denote the space $\ell_\infty(\Omega)$ quasi-normed
by $\|\cdot\|_p$.
\begin{lemma}\label{kalton}
\begin{itemize}
\item[(c)] $X_p(\Omega)$ is a $p$-normed lattice.
\item[(d)] $\|f\|_\infty\leq\|f\|_p$ for every $f\in\ell_\infty(\Omega)$.
\item[(e)] $\|1_\Omega\|_p\geq (n^{1/p-1})\varepsilon^{1/p}/(1+\varepsilon)$.
\item[(f)]
$\|f\|_\infty\leq\|f\|_{\co(X_p(\Omega))}\leq(1+\varepsilon)\|f\|_\infty$
for all $f$.
\end{itemize}
\end{lemma}
\begin{proof} (c) and (d) are trivial. To verify (e), let $c_i\geq0$ be
so that $1_\Omega\leq\sum_{i=1}^{(1+\epsilon)n}c_i1_{A_i}$. Put
$$
J=\left\{i: c_i\geq\frac{1}{(1+\varepsilon)n}\right\}.
$$
We claim that $|J|>\varepsilon n$. Suppose on the contrary that
$|J|\leq\varepsilon n$. By Lemma \ref{talagrand}, there is
$\omega\in\Omega$ so that $\omega\notin\cup_{i\in J}A_i$. Hence
$$
\sum_{i=1}^{(1+\varepsilon)n}1_{A_i}(\omega)=\sum_{i\notin
J}1_{A_i}(\omega)<\sum_{i\notin J}\frac{1}{(1+\varepsilon)n}<1,
$$
a contradiction. Therefore $|J|\geq\varepsilon n$ and
$$
\sum_i c_i^p\geq\varepsilon n \left(\frac{1}{(1+\varepsilon)n}\right)^p
=\frac{\varepsilon}{(1+\varepsilon)^p}n^{1-p},
$$which proves (e).

The first inequality in (f) is clear. As for the other, note that
$\co X_p(\Omega)$ is a Banach lattice in its natural order. Let
$f$ be such that $\|f\|_\infty\leq1$. Then $|f|\leq1_\Omega$ and
since $1_\Omega=\frac{1}{n}\sum_{i=1}^{(1+\varepsilon)n}1_{A_i}$
one has
$$
\|f\|_{\co X_p(\Omega)} \leq \|1_\Omega\|_{\co X_p(\Omega)}
\leq\frac{1}{n}\sum_{i=1}^{(1+\varepsilon)n}\|1_{A_i}\|_{X_p(\Omega)}
=1+\varepsilon,
$$ from which the result follows.
\end{proof}

To end with the example, let $\varepsilon>0$ be any rational number. Choose
a sequence $(n_k)_k$ such that $\varepsilon n_k\in\mathbb N$ for all
$k\in\mathbb N$. Fix $0<p<1$ and define $X$ to be the $c_0$-sum of the
spaces $X_p(\Omega(\varepsilon,n_k))$, that is,
$$
X=c_0(X_p(\Omega(\varepsilon,n_k))=\left\{(f_k)_k: f_k\in
 X_p(\Omega(\varepsilon,n_k)),\text{ with }\lim_k\|f_k\|_p=0\right\}
$$
equipped with the $p$-norm
$$
\|(f_k)\|_ X= \sup_k \|f_k\|_{p}.
$$
Clearly, $X^*$ separates points in $X$. Since
$$
\co X=\co(c_0(X_p(\Omega(\varepsilon,n_k)))=c_0(\co
(X_p(\Omega(\varepsilon,n_k)))
$$
it follows from Lemma \ref{kalton} that $X$ is not locally convex.
On the other hand (f) shows that the convex envelope is
$(1+\e)$-isomorphic to $c_0$.
\end{proof}

\end{document}